\documentclass{article}
\usepackage[dvipdfmx]{graphicx}
\usepackage{amsmath, amssymb}
\usepackage{type1cm}
\usepackage{bm}
\usepackage{here}
\usepackage{mathrsfs}
\usepackage[margin=2cm]{geometry}
\usepackage{wrapfig}
\usepackage{authblk}

\newtheorem{theorem}{Theorem}

\newtheorem{ho}{Lemma}

\title{Convergence of Q-value in case of Gaussian rewards}
\makeatletter
  
  \@addtoreset{equation}{section}
\makeatother
\author[1,3]{Konatsu Miyamoto}
\author[4,1]{Masaya Suzuki}
\author[2]{Yuma Kigami}
\author[3,1]{Kodai Satake}

\affil[1]{Dynamic Pricing technology}
\affil[2]{MatrixFlow}
\affil[3]{Osaka univercity}
\affil[4]{Kindai univercity}

\begin{document}
\maketitle

\begin{abstract}
In this paper, as a study of reinforcement learning, we converge the Q function to unbounded rewards such as Gaussian distribution.
From the central limit theorem, in some real-world applications it is natural to assume that rewards follow a Gaussian distribution , but existing proofs cannot guarantee convergence of the Q-function. Furthermore, in the distribution-type reinforcement learning and Bayesian reinforcement learning that have become popular in recent years, it is better to allow the reward to have a Gaussian distribution. 
Therefore, in this paper, we prove the convergence of the Q-function under the condition of $E[r(s,a)^2]<\infty$, which is much more relaxed than the existing research.
Finally, as a bonus, a proof of the policy gradient theorem for distributed reinforcement learning is also posted. 
\end{abstract}

\section{Introduction}

In recent years, Reinforcement Learning(RL) has come into fasion. 
General method in ordinary Reinforcement Learning using Markov decision processes use a state action value functions[1].
Agents created by these algorithms take strategies to maximize the expected value of the cumulative reward.
However, in practical use , there are many situations where it is necessary to consider not only expected values but also risks.
Therefore, Distributional Reinforcement Learning(DRL) that considers the distribution of cumulative rewards has also been studied.
DRL research presents a particle method of risk responsive algorithm[2].
As for similar research, there are[3][4],which is equivalent to [2] mathematically,but used the different algorithm and parametric methods[5]. 
[4]  discusses the convergence of measures in discrete steps.
Another way to practice DRL is using the Bayesian approach.
In [22],it is regarded as an estimation of the uncertainty of the expected value.But in fact, the Bayesian inferece can approximate the distribution of uncertain objecsion.
It can perform distributed reinforcement learning.
There are other  existing papers on Bayesian reinforcement learning.
We want to take [6][7] up this time.It is a method using Gaussian processes, and it can be said that the reward follows Gaussian distributions.
[5] also supports unbounded rewards like Gaussian distributions.
We want to show that the approximation of the cumulative reward distribution converges even in unbounded rewards.
In this paper, we prove the convergence of the normal state action value function as a preliminary step.
In addition, we perform the convergence proof for Q functions with continuous concentration domain,taking Deep Q-learning(DQN) into consideration.

\subsection{Related works}
The proof history of Q-function convergence is long.
For example, there are papers such as [8], [9], [10], and [11] using [10].
A paper on an unusual proof method is [12] using ordinary differential equations.
For DQN, there is a study [13] summarizing the approximation error.
The approximation error due to the neural network is verified there.
Other research results include [14][15][16][17][18].
All of these studies assume that rewards are bounded. That is, there is a certain constant $R_{max}<\infty$ and

\begin{align}
|r(s,a)|\leq R_{max}\ a.e.
\end{align}
holds.
Therefore, Gaussian distributions cannot be assumed.
In this paper, we prove the convergence of the Q function under condition ,
\begin{align}
\forall s,a\in S\times A,E[r(s,a)^2]<\infty
\end{align}
which is more relaxed than (1.1), with normal distribution in mind.
Finally, we prove the convergence of the Q function in the domain of continuous concentration under ideal conditions.
This is a frequent concept in reinforcement learning.

\section{Background}

\subsection{transition kearnel}

Let two tuples $(S,\mathcal{S}),(T,\mathcal{T}) $ be both measurable spaces.Transition kernel $k:S\times \mathcal{T}\to\mathbb{R}_+$ is defined to satisfy the following two conditions.
\begin{align}
\cdot& \forall B\in\mathcal{T}, k(\cdot,B)\ on\ S\ is\ measurable\\
\cdot& \forall s\in S,k(s,\cdot)is\ measure\ on\ T
\end{align}
This is used in situations where $ s $ is fixed and the distribution on $ T $ is fixed.

\subsection{Markov decision process}
Assume that both the set of states $S$ and the set of actions $A$ are finite sets.
A transition kernel $p$ is defined on $(S\times A,2^{S\times A}),(S\times \mathbb{R},2^S\otimes\mathcal{B}(\mathbb{R}))$.
That is, $p(r,s|s,a)$ is a probability measure that governs the distribution of the next state $s\in S$ and immediate reward $r \in \mathbb{R}$ when an action $a\in A$ is taken in state $ s \in S$ is there.
The strategy $ \pi: S \to \mathcal{P}(A) $ is the action probability determined from the current situation, as can be seen from the definition. The deterministic approach is that for any $ s$, there is a $ a $ and $ \pi (a | s) = 1 $.
A set of random variables $ s_t, a_t, r_t $ taking values in $ S, A, \mathbb{R} $ is written as $ (s_t, a_t, r_t)^\infty_{t = 0} $.
This stochastic process is called Markov decision process(MDP).
\subsection{Optimal measures and state action value functions}
Put the whole set of policies as $ \Pi $.
The state action value function $ Q^\pi: S \times A \to \mathbb{R} $ for the policy $ \pi $ is defined as follows.
\begin{align}
Q^\pi(s,a):=E[\sum^\infty_{t=0}\gamma^t R_t|s_0=s,a_t=a(r_t,s_{t+1}),p(r_t,s_{t+1}|s_t,a_t),\pi(a_t|s_t)]
\end{align}
Furthermore, the state value function $ V^\pi (s) $ is defined as follows.

\begin{align}
V^\pi(s):=\sum_{a\in A} \pi(a|s)Q^\pi(s,a)
\end{align}
Define the optimal strategy $ \pi^* $ as
\begin{align}
\pi^*:=argmax_{\pi\in\Pi} V^\pi(s_0) 
\end{align}
In addition, the state action value function $ Q^{\pi^*} $ for the optimal policy is called the optimum state action value function, and simply expressed as $ Q^* $.
The action that takes the maximum value for the optimal state action function is the optimal policy. 

\begin{align}
\pi^*(a|s)=\begin{cases}
1\ argmax_{a\in A} Q(s,a)\\
0\ else
\end{cases}
\end{align}
holds for any $ s, a $.

\section{Update of state action value function and Robbins Monro condition}

Update the Q-unction as follows
\begin{align}
Q_{n+1}(s,a)=(1-\alpha(s,a,s_t,a_t,t))Q_t(s,a)+\alpha(s,a,s_t,a_t,t)[r_t(s_t,a_t)+\max_{b\in A}Q_t(s_{t+1},b)]
\end{align}
The following sequence $ \{c_t \}^\infty_{t = 0} $ satisfies the Robbins-Monro condition.

\begin{align}
\forall t, c_t\in[0,1]\\
\sum^\infty_{t=0} c_t=\infty\\
\sum^\infty_{t=0}c^2_t<\infty
\end{align}
Using this, the mapping $ \alpha: S \times A \times S \times A \times \mathbb{N} \to (0,1) $ is defined as follows.

\begin{align}
\alpha(s,a,s_t,a_t,t)=\begin{cases}
c_t \ s_t=s,a_t=a\\
0\ else
\end{cases}
\end{align}
In addition,it is assumed that this also satisfies the Robbins Monroe condition stochastically uniformly for arbitrary $ s, a $.
\begin{align}
\sum^\infty_{t=0} \alpha(s,a,s_t,a_t,t)=\infty\ a.e.\\
\sum^\infty_{t=0}\alpha(s,a,s_t,a_t,t)^2<\infty\ a.e.
\end{align}

\section{Proof of Q-function convergence for unbounded rewards}

Consider a real-valued function $ w_t (x) $ on a finite set $ \mathcal{X} $.
\begin{theorem} Convergence of Q-value in case of Gaussian rewards

$\mathcal{X}$ is finite set.
Let ramdom value $r_t(x), \mathcal{X}: = S \times A $.
Let $ W $ be the set of functions $ f: \mathcal{X} \to \mathbb{R} $ and $||f||_W$is defined as$ || f ||_W: = \max_ {x \in \mathcal{X}} f (x)$ .
For any $ s, a $, let $ E [r^2 (s, a)] <\infty $.

\begin{align}
||Q_t-Q^*||_W \to 0
\end{align}
proof.

In line with the proof of [9]. The $ F $ condition is relaxed and the statement is stronger, so it needs to be done more precisely.
Consider a stochastic process of $ \Delta_t (x): = Q_t (x) -Q^* (x) $.
Since $ Q^* (x) $ is a constant, $ V (\Delta_t (x)) = V (Q_t (x)) $.
Putting$ F_t (x): = r_t (x) + \gamma \sup_b Q_t (X (s, a), b) -Q^* (x) $, this is $ \mathcal{F}_{t + 1} $ measurable stochastic process.
Furthermore, if we put $ G_t (x): = r_t (x) + \gamma \sup_b Q_t (X (s, a), b) $, by definition $ G_t-E [G_t (x) | \mathcal { F}_t] = F_t-E [F_t (x) | \mathcal{F}_t] $.
The two stochastic processes $ \delta_t, w_t \in W $ are taken so that $ \Delta_0 (x) = \delta_0 (x) + w_0 (x) $.
Define time evolution as
\begin{align}
 &\delta_{t+1}(x)=(1-a_t(x))\delta_t(x)+a_t(x)E[F_t(x)|\mathcal{F}_t]\\
& w_{t+1}(x)=(1-a_t(x))w_t(x)+a_t(x)p_t(x)
\end{align}
However, $ p_t (x): = F_t (x) -E [F_t (x) | \mathcal{F}_t] $. At this time, $ \Delta_t (x) = w_t (x) + \delta_t (x) $.
First, we show that $ w_t $ converges to 0 for $ \mathcal{X} $ with probability 1 by using Lemma 2.
By definition, $ E [p_t | \mathcal {F}_t] = 0 $, so $ \sum_t E |[p_t | \mathcal{F} _t]| = 0 $ holds.
From Lemma 1 and the definition of $ p_t, G_t $, $ E [p_t^2] \leq 4E [G_t^2]$ holds.
Putting $ L_t (\omega): = \sup_x | Q_t (x) | $, this random variable is $ \mathcal{F}_t $ -measurable and takes a finite value with probability $ 1 $.
Since $ L_0 $ is a finite value, a certain constant $ K_0 $ can be taken so that $ E [L_0^2] \leq K^2_0C_R $ holds.
And the following holds with probability 1.

\begin{align}
L_{t+1}\leq \max(L_t,(1-b_t)L_t+b_t(\sup_x|r_t(x)|+\gamma L_t))
\end{align}
Using the above formula, the following holds
\begin{align}
E[L^2_{t+1}]&\leq \max(E[L^2_t],E[((1-b_t)L_t+b_t(\sup_x|r_t(x)|+\gamma L_t))^2])
\end{align}

Suppose there is $ K_t \in \mathbb{R} $ that is $ E [L_t ^ 2] \leq K^2_tC_R $. At this time, put $ H_t: = \sup_x | r_t (x) | + \gamma L_t $

\begin{align}
E[H_t^2]&=E[\sup_x|r_t(x)|^2]+2E[\sup_x|r_t(x)|\gamma L_t]+\gamma^2E[L_t^2]\\
&\leq C_R+2\gamma\sqrt{C_RK_t^2C_R}+K^2_tC_R\\
&=(1+\gamma K_t)^2C_R
\end{align}
Then,
\begin{align}
E[((1-b_t)L_t+b_t(\sup_x|r_t(x)|+\gamma L_t))^2]&\leq (1-b_t)^2E[L_t^2]+2(1-b_t)b_t\sqrt{E[L_t^2]E[H_t^2]}+b_t^2\gamma^2E[H_t^2]\\
&\leq (1-b_t)^2K_t^2C_R+2(1-b_t)b_tK_t(1+\gamma K_t)C_R+(1+\gamma K_t)^2C_R\\
&=((1-b_t)K_t+b_t(1+\gamma K_t))^2C_R\\
&=(K_t+b_t(1-(1-\gamma)K_t))^2 C_R
\end{align}
Putting $ K_{t + 1} = max (K_t, K_t + b_t (1- (1- \gamma) K_t)) $, $ E {L^2_{t + 1}} \leq K_{t + 1} C_R $ can be said. Since $ K_0 \in \mathbb {R} $ exists, $ K_t \in \mathbb {R} $ exists for any $ t $, and $ E [L_t^2] \leq K_t^2 C_R$ can be said.
It is clear from the equation that $ K_{t + 1} = K_t $ when $ K_t> \frac {1} {1- \gamma} $, and $ K_t \leq \frac {1} {1- \gamma} $ Then, $ K_{t + 1} \leq \frac {1} {1- \gamma} + 1 $ holds.
Therefore, it was shown earlier that $ K_t $ exists for any $ t $, in addition $ K_t \leq K^*: = \max (K_0, \frac {1} {1- \gamma} + 1) $ can be also said.
$ | G_t (x) | \leq | r_t (x) | + \gamma L_t $ holds, so the following equation hold.for all $x$
\begin{align}
\frac{1}{4}E[p_t^2(x)]&\leq E[G_t^2(x)]\\
&=E[r_t(x)^2]+2\gamma \sqrt{E[r_t(x)^2]E[L_t^2]}+E[L_t^2]\\
&\leq (1+\gamma K^*)C_R
\end{align}
Then,
\begin{align}
\sum_tE[a_t^2p_t^2]&\leq \sum_t 4b^2_t (1+\gamma K^*)C_R\\
&\leq4 M (1+\gamma K^*)C_R<\infty
\end{align}
holds for all $x$.
When we use Lemma2,putting
\begin{align}
&U_t:=a_t(x)p_t(x)\\
&T(w_t,\omega):=(1-a_t(x))w_n
\end{align}
$ \sum_t E [U_t^2] <\infty$can be said.
Since $ E [U_t | \mathcal{F}_n] = 0 $, $ \sum_t | E [U_t | \mathcal {F}_n] | = 0 $ holds.
Then, for any $ \epsilon> 0 $, set $ \alpha = \epsilon, \beta_t (x) = b^2_t (x) $ and $ \gamma_t (x) = \epsilon (2a_t ( x) -a_t^2 (x)) $, then 
\begin{align}
&T^2(w_t,\omega)\leq \max(\alpha,(1+\beta_t)w^2_t-\gamma_t)\\
&\sum_t\gamma_t =\infty\ a.e 
\end{align}
holds.
The latter is based on Robbins Monro conditions.
Therefore, $ w_t (x) \to0 $ holds for any $ x $.
Define the linear operator $ \mathcal{T}: W \to W $ as follows: for $ q \in W $
\begin{align}
\mathcal{T}q(s,a)&=\int_\mathbb{R}\sum_{s'} [r(s,a)+\gamma \sup_b q(s',b) ]p(dr,s'|s,a)\\
&=E[r(s,a)+\sup_bq(X(s,a),b)]
\end{align}
$ Q^* $ is a fixed point for this operator.
For any $q_1,q_2\in W$
\begin{align}
||\mathcal{T}q_1-\mathcal{T}q_2||_W&=sup_{s,a}[|\int_\mathbb{R}\sum_{s'} [r(s,a)+\gamma \sup_b q_1(s',b) ]p(dr,s'|s,a)-\int_\mathbb{R}\sum_{s'} [r(s,a)+\gamma \sup_b q_2(s^`,b) ]p(dr,s'|s,a)|]\\
&\leq\int_\mathbb{R}\sum_{s'} [\gamma |\sup_b q_1(s',b)-\sup_bq_2(s',b)| ]p(dr,s'|s,a)\\
&\leq \int_\mathbb{R}\sum_{s'} [\gamma sup_b |q_1(s',b)-q_2(s,b)| ]p(dr,s'|s,a)\\
&=\gamma ||q_1-q_2||_W
\end{align}
Thus $ \mathcal{T} $ is a reduction operator.
\begin{align}
|E[F_t(x,a)|\mathcal{F}_t]|&\leq\int_\mathbb{R}\sum_{s'} |r(s,a)+\gamma \sup_b Q_t(s',b) -Q^*(s,a)|p(dr,s'|s,a)\\
&=|\mathcal{T}Q_t(x,a)-Q^*(s,a)|\\
&=|\mathcal{T}Q_t(x,a)-\mathcal{T}Q^*(s,a)|\\
&\leq \gamma ||\Delta_t||_W
\end{align}
Then,
\begin{align}
||\delta_{t+1}||&\leq (1-a_t(x))||\delta_t||+a_t(x)||\delta_t+w_t||\\
&\leq(1-a_t(x))||\delta_t||+a_t(x)(||\delta_t||+||w_t||)
\end{align}
$ || w_t (x) || $ converges uniformly to 0 with a probability of 1 for any $ x $ as described above.
Therefore, from Lemma 3, $ || \delta_{t + 1} (x) || \to0 $ for any $ x $.
That is, for any $ x $, $ || \Delta_t (x) ||_W \to0 $, which holds the main theorem assertion.

\end{theorem}

\section{Theorem for SARASA}

The method in Chapter 3 is called Q-learning, and the value is updated before performing the next action.
On the other hand, SARASA updates the value after performing the following actions.
\begin{align}
Q_{t+1}(s,a)=(1-\alpha(s,a,s_t,a_t,t))Q_t(s,a)+\alpha(s,a,s_t,a_t,t)(r_t(s,a)+Q_t(s_{t+1},a_{t+1}))
\end{align}
$ a_{t + 1} $ is often stochastically determined by softmax function or the like.
\begin{theorem}
Suppose that the Q function is updated by the above SARASA method. At this time,
\begin{align}
||Q_t-Q^*||_W\to0\ in \ t\to\infty
\end{align}

proof.

Put $ L'_t: = \ max_ {x, y \ in \ mathcal {X}} | Q_t (x) -Q_t (y) | $ It is clear from the definition that $ L'_t \leq 2L_t $.
Later along this follows the proof of Theorem 1.

\section{Convergence proof for unbounded rewards under continuous concentration}
For example, in a situation such as DQN, an update for one $ s, a $ has an effect on other state actions.
As a simple model to take such situations into account, we put the ripple function $ f (x_1, x_2) $ defined on the compact set $ \mathcal {X}^2 $. This satisfies the next conditions.

\begin{align}
&f(x,x)=1\\
&f(x,y)\ is\ continue.
\end{align}
\end{theorem}
If $ Q^* $ is a continuous function, it can be used to depart from any continuous function and have the same convergence on the compact set.
Let $ \mathcal {X} \subset \mathbb{R}^d $ be a simple connected compact set.
Let $ Q^*, Q_0 $ be a continuous function on $ \mathcal{X} $.
Let $ W $ be a continuous function on $ \mathcal{X} $. $ || f ||_W: = \max_{x \in \mathcal {X}} f (x) $
\begin {align}
Q_{t + 1} (s, a) = (1-f (s, a, s_t, a_t) \alpha (s, a, s_t, a_t, t)) Q_t (s, a) + f (s, a, s_t, a_t) \alpha (s, a, s_t, a_t, t) (r_t (s, a) + \max_ {b \in A} Q_t (s_{t + 1}, b))
\end {align}
At this time, $ || Q_t-Q^* ||_W \to 0 $

proof.
Consider a finite set $ K_N: = \{x_1, x_2, x_3, ......, x_N \} $ on $ \mathcal{X} $. Limiting $ Q $ to $ K $ converges to a correct function uniformly over $ K $ from Theorem 1.For any $ \epsilon $
Since $ Q^* $ is a continuous function, the function whose value is defined on a dense set is uniquely determined.
Convergence can be said.

\section{Conclusion and Future Work}

As we mentioned earlier,we want to prove the convergence of the distribution. An order evaluation of the expected value should be also  performed.
We also want to estimate the convergence order for a specific neural network such as [13]. According to [13], as with Theorem 3, in the domain of continuous concentration, as $ R_ {max}: = \sup r(\omega, s, a) $, using constants $ C_1, C_2, \xi, \alpha $ 
\begin{align}
||Q^*-Q_n||_W\leq C_1\cdot (\log n)^\xi n^{-\alpha}+C_2R_{max}
\end{align}
is established.
However, when $ r $ follows a normal distribution, $ R_{max} = \infty $, so the upper limit of the error is infinite, and this unexpected expression has no meaning. In case of using unbounded rewards, stronger inequality proofs are needed.

\appendix \def\thesection{\Alph{section}}

\section{Lemmas and proofs}
\begin{ho}

　

Consider a random variable $ Y $ and a partial $\sigma$-algebla $ \mathcal{G} $.
If $ Z: = Y-E [Y | \mathcal{G}] $, the following equation holds.

\begin{align}
E[Z^2]\leq 4E[Y^2]
\end{align}

\end{ho}

We quote the important theorem.

\begin{ho}Convergence theorem for stochastic systems[19]

　

Consider the following stochastic process.
\begin{align}
X_{t+1}:=T(X_0,......,X_t,\omega)+U_t(\omega)
\end{align}
This satisfies the following equation with probability 1.

\begin{align}
|T(x_1,x_2,......,x_t,\omega)|^2\leq max(\alpha,(1+\beta_t(\omega))x^2_t-\gamma_t)
\end{align}
However, with $ \alpha> 0 $, with probability 1, $ \beta_t (\omega) <M ', \sum_t \beta_t <\infty $ holds, and with probability 1, Let $ \sum_t \gamma_t (\omega) = \infty $.
\begin{align}
\sum_t E[U_t^2]&<\infty\\
\sum_t E[U_t|\mathcal{F}_t]&<\infty
\end{align}
At this time, there exists a certain $ N (\omega) $, and it  holds for any $ n> N (\omega) $
\begin{align}
\lim\sup_{t\to \infty}|X_t|^2<\alpha\ a.e.
\end{align}

\end{ho}
If $ \beta,\gamma $ are taken again for any $ \alpha $ and the same can be said, "uniform convergence to 0" can be said that is much stronger than approximate convergence.

\begin{ho}

$ x_0\in \mathbb{R} $ is assumed to be a real number.

\begin{align}
x_{n+1}=(1-a_n)x_n+\gamma a_n|x_n|
\end{align}

$ \gamma \in (0,1) $ is a constant.At this time, $ x_n \to0 $ holds with probability 1.

proof.

Look at each $ \omega $. That is, $ \{a_n \}^\infty_{n = 0} $ is constant sequence that satisfies $ \sum^\infty_{n = 0} a_n = \infty, \sum^\infty_{n = 0} a^2_n <\infty $.
$ X_n $ is nonnegative for a sufficiently large $ n $, so it is bounded below.
In addition, since $ x_n \geq x_{n + 1} $ is apparent from the equation, $ \{x_n \}_{n = 1}^\infty $ is a monotonically decreasing sequence.
The sequence converges because it is bounded and monotonically decreasing below.Putting
$ b_n: = a_n- \gamma a_n $, this satisfies $ \sum^\infty_{n = 0} b_n = \infty, \sum^\infty_{n = 0} b^2_n <\infty $ .
You can say$ x_n = x_1 \Pi^n_{i = 1} (1-b_i) $, and the convergence destination $ x $ is $ x_1 \Pi^\infty_{n = 1} (1-b_n) $.
$ c_n = \Pi^n_{i = 1} (1-b_i) $, the infinite product of $ n \to \infty $ is $ \sum^\infty_{n = 0} b_n = \infty$, but diverges. However, since it is $ 0 \leq c_n \leq 1 $, $ c_n \to0 $ is known, and $ x_n \to0 $ can be said.
\end{ho}

\begin{ho}

Let $\epsilon>0$.

\begin{align}
x_{n+1}=(1-a_n)x_n+\gamma a_n|x_n+\epsilon|
\end{align}
Then $x_n\to \epsilon\frac{\gamma}{1-\gamma}$ holds.

proof.

\begin{align}
x_{n+1}-x_n&=-a_n((1-\gamma)x_n-\epsilon \gamma) \\
&=-a_n(1-\gamma)(x_n-\epsilon\frac{\gamma}{1-\gamma})
\end{align}
The difference from $ \epsilon \frac {\gamma} {1- \gamma} $ is reduced by $ a_n (1- \gamma) $.
If $ y_n: = x_n- \epsilon \frac {\gamma} {1- \gamma} $, by definition it is clearly $ y_{n + 1} -y_n = x_{n + 1} -x_n $.
Moreover,
\begin {align}
y_{n + 1} -y_n & =-a_n (1- \gamma) (y_n) \\
y_{n + 1} & = (1-a_n (1- \gamma)) y_n
\end {align}
After that, it is $ x_n \to \epsilon \frac{\gamma}{1- \gamma} $ because it is $ y_n \to0 $ by the same argument in Lemma 3.
\end{ho}
\begin{ho}
Suppose that the sequence $ \{c_n\} \subset \mathbb{R}_+ $ converges uniformly to 0 on a set of probabilities 1.
That is, for any $ \epsilon_1> 0 $, there is a certain $ N_{\epsilon_1}(\omega) $, and when $ n> N_{\epsilon_1} (\omega) $, $ | c_n | <\epsilon_1 $ holds with probability 1.
At this time,
\begin {align}
x_ {n + 1} = (1-a_n) x_n + \gamma a_n | x_n + c_n |
\end {align}
${x_n}$ converges to 0.

proof.

\begin{align}
z_{N_{\epsilon_1}}&=x_{N_{\epsilon_1}}\\
z_{n+1}&=(1-a_n)z_n+\gamma a_n|z_n+\epsilon_1|
\end{align}
$ | Z_n | \geq | x_n | $ for such $ n> N $. $ Z_n \to \epsilon_1 \frac{\gamma} {1- \gamma} $ from Lemma 4.
That is, for any $ \epsilon_2> 0 $, there is a certain $ N_{\epsilon_2}> N_{\ epsilon_1} $, and $ n> N_{\epsilon_2} $ for any $ n , z_n <\epsilon_1 \frac{\gamma}{1- \gamma} + \epsilon_2 $
$ \epsilon_1, \epsilon_2 $ can be arbitrarily taken, so if we define a new $ \epsilon: = \epsilon_1 \frac{\gamma}{1- \gamma} + \epsilon_2 $, this is also $ \epsilon> 0 $ can be taken arbitrarily. Within the range of Using $ z_n> x_n $, there is a $ N_{\epsilon_2} $ for any $ \epsilon $ and $ x_n <\epsilon $for $ n> N_{\epsilon_2} $ .

\end{ho}

\section{Strict Proof of Policy Gradient thorem and Distributionaly}

We prove the famous policy gradient theorem using the $ Q $ function and its version in distributed reinforcement learning [23].

\begin{theorem}Policy Gradient thorem

　

Consider the gradient of the policy value function $ J (\theta):=E [Q (x, \pi_\theta (x))] $. At this time, it is assumed that $ \pi, Q $ is implemented by a neural network, the activation function is Lipschitz continuous, and $ \nabla_\theta Q (x, a) = 0 $.
Then,The following equation holds,
\begin{align}
\nabla_\theta J(\theta)=E_\rho[\nabla_\theta \pi(\theta) \nabla_a Q(s,a)|_{a=\pi(x)}]
\end{align}
However, $ \rho $ is memory data in general implementation.
Next, consider the case of distributed reinforcement learning. If a random variable representing the cumulative reward sum is expressed as $ Z $, then $ Q (s, a) = E_n [Z (s, a)] $ holds. Suppose $ Z $ is a neural network with stochastic output. 

\begin {align}
Z (\omega) (s, a) = f_\omega (s, a)
\end{align}

Then
\begin{align}
\nabla_\theta J(\theta)=E_\rho[\nabla_\theta \pi(x) E_n[\nabla_a Z(x,a) ]|_{a=\pi(x)}]
\end{align}

proof.

The interchangeable conditions of differentiation and Lebesgue integration are described as follows.
Suppose there is a function $ f (x, \omega) $ that can be Lebesgue integrable over $ \Omega $ and differentiable by $ x $.
At this time, there is an integrable function $ \phi (\omega) $, and $ x $ can be differentiated almost everywhere on $ \Omega $ by $ x $ and $ | \nabla_x f (x, \omega) _i | \leq \phi (\omega) $ holds, then $ \int_\Omega f (x, \omega) d \mu (\omega) $ is differentiable by $ x $, and holds,
\begin{align}
\nabla_x \int_\Omega f(x,\omega)d\mu(\omega)=\int_\Omega \nabla_xf(x,\omega)d\mu(\omega)
\end{align}
When $\mu(\Omega)<\infty$, An example of a function class that satisfies this is the “Lipschitz continuous function”. Neural networks is generally combinations of linear transformations and Lipschitz continuous activation map.\footnote{
General activation functions such as sigmoid, ReLu, Reaky Relu, and Swish are all Lipschitz continuous functions.}
Moreover, if the Lipschitz constant of the function $ f $ is written as $ || f ||_L $, then considering two Lipschitz continuous functions $ f, g $, $ || f \circ g ||_L \leq | | f ||_L || g ||_L $.
From this, $ \pi_\theta (x), Q (x, a), Q (x, \pi_\theta (x)) $ are Lipschitz continuous for $ x, a $, respectively. Although it is not Lipschitz continuous for $ \theta $, it is Lipschitz continuous for each element, and the definition and definition of $ \nabla $ allow the exchange of differentiation and integration.
That is, the following holds from the differential chain rule,

\begin{align}
\nabla_\theta J(\theta)=E_\rho[\nabla_\theta \pi_\theta(x) \nabla_a Q(s,a)|_{a=\pi_\theta(x)}]
\end{align}
Similarly, $ \nabla_a E [Z (s, a)] = E [\nabla_a f_\omega (s, a)] $ and $ f_\omega $ is Lipschitz continuous functions for any $ \omega $ , For distribution type
\begin{align}
\nabla_\theta J(\theta)&=E_\rho[\nabla_\theta \pi(x) E_n[\nabla_a f_\omega (x,a) ]|_{a=\pi(x)}]\\
&=E_\rho[\nabla_\theta \pi(x) E_n[\nabla_a Z (x,a) ]|_{a=\pi(x)}]
\end{align}
\end{theorem}
As described above, the policy gradient theorem is established because the policy is Lipschitz-continuous for each parameter, and is obviously not for a policy function composed of ODEnet[24], hypernet[25], or the like that reuses parameters.

\section{Notation}
Let $ (A, \mathcal{O}) $ be a topological space.
\begin{itemize}
\item $\sigma(\mathcal{O})$:Smallest $\sigma$-algebla containing all $o\in\mathcal{O}$.
\item $\mathcal{B}(A):=\sigma(\mathcal{O})$
\item $\mathcal{P}(A)$:If $ A $ is a finite set, it is the set of all probability measures defined by measurable space $ (A, 2^A) $, and if it is an infinite set, $ (A, \mathcal{B} (A)) $
\item $\mathcal{F}\otimes\mathcal{G}:=\sigma(\mathcal{F}\times\mathcal{G})$
\end{itemize}

\end{document}